\documentclass[a4paper]{amsart}

\usepackage{amsthm}
\usepackage[T1]{fontenc}
\usepackage[utf8]{inputenc}
\usepackage[hidelinks]{hyperref}
\usepackage{needspace}

\newtheorem{Theorem}{Theorem}
\newtheorem{Corollary}[Theorem]{Corollary}
\newtheorem{Lemma}[Theorem]{Lemma}

\theoremstyle{definition}

\newcommand{\C}{\mathbb{C}}
\newcommand{\Z}{\mathbb{Z}}
\newcommand{\N}{\mathbb{N}}
\newcommand{\R}{\mathbb{R}}
\newcommand{\floor}[1]{\lfloor #1  \rfloor}

\begin{document}

\title[Second moments of general Dirichlet series]
{Second moments of general Dirichlet series and Selberg zeta functions at $\sigma = 1$}

\author[R. Garunkštis]{Ramūnas Garunkštis}
\address{Ramūnas Garunkštis\\ Institute of Mathematics\\
Faculty of Mathematics and Informatics\\
Vilnius University\\
Naugarduko 24, 03225 Vilnius, Lithuania}
\email{ramunas.garunkstis@mif.vu.lt}

\author[J. Putrius]{Jok\={u}bas Putrius}
\address{Jok\={u}bas Putrius \\
 Institute of Mathematics, Faculty of Mathematics and Informatics, Vilnius University \\
Naugarduko 24, 03225 Vilnius, Lithuania}
\email{jokubas.putrius@mif.vu.lt}

\keywords{General Dirichlet series, second moment, Selberg zeta-function, Beurling zeta-function}
\subjclass[2020]{11M41}

\begin{abstract}
In this paper, we prove the existence of the second moment at $\sigma=1$ for a general Dirichlet series with positive coefficients whose partial sum of the coefficients up to $x$ satisfies
$\varrho x+O(x/(\log x)^{1+\varepsilon})$ for some $\varepsilon>0$.
The result extends a theorem of Broucke and Hilberdink.
As applications, we establish the existence of the second moment for Beurling zeta functions satisfying a prime number theorem with a power-saving error term and for the Selberg zeta function for every cofinite Fuchsian group.
\end{abstract}

\maketitle

\section{Introduction}

As usual, we write $s=\sigma+it$. In this paper, we investigate the second moment at $\sigma=1$ of a general Dirichlet series with positive coefficients. Our main result is the following theorem.

\begin{Theorem}\label{thm:moment2}
Let
\begin{equation}\label{f dirichlet series}
f(s)=\sum_{j=1}^{\infty}\frac{a_j}{n_j^s},
\end{equation}
where $(n_j)$ is a strictly increasing sequence of positive reals tending to infinity and $a_j>0$. Assume that
\begin{equation}\label{eq:ax}
A(x)=\sum_{n_j\le x}a_j=\varrho x+R(x),
\end{equation}
where $\varrho\ge0$ is a constant and, for some fixed $\varepsilon>0$,
\begin{equation}\label{eq:Rxbound}
R(x)\ll\frac{x}{(\log x)^{1+\varepsilon}}.
\end{equation}
Then the Dirichlet series defining $f(s)$ converges absolutely for $\sigma>1$, and $f(s)$ has a continuous extension to $\sigma\ge1$, $s\ne1$. Moreover,
\begin{equation}\label{mainmoment}
\lim_{T\to\infty}\frac1T\int_1^T|f(1+it)|^2\,dt
=\sum_{j=1}^{\infty}\frac{a_j^2}{n_j^2}.
\end{equation}
The series on the right-hand side converges.
\end{Theorem}

Mean values of general Dirichlet series were studied by Landau \cite[Chapter 75]{landau1909}, whose results require a separation condition on the exponents. Recently, Broucke and Hilberdink \cite[Theorem 1.1]{bh} proved a second moment theorem using an alternative condition. They assume that the coefficients $a_j$ are positive and satisfy
\[
\sum_{n_j\le x}a_j=\varrho x+O(x^\beta),\qquad \varrho>0,\quad 0\le\beta<1.
\]
They then establish the convergence of the second moment in $\sigma>(1+\beta)/2$. Theorem \ref{thm:moment2} treats the case $\sigma=1$ under a weaker remainder estimate, which does not provide a power saving error term.

The following theorem was inspired by a question of F. Broucke during the ELAZ2026 conference. It shows that Theorem \ref{thm:moment2} is sharp in the sense that replacing $1+\varepsilon$ by $1$ in \eqref{eq:Rxbound} makes the conclusion false in general. Let $\{x\}$ denote the fractional part of $x$.

\begin{Theorem}\label{claim}
Let $f(s)=\sum_{j=1}^{\infty}a_jn_j^{-s}$ with
\[
n_j=e^{\sqrt{j}}\qquad\text{and}\qquad
a_j=\frac{e^{\sqrt{j}}}{\sqrt{j}}.
\]
Then
\begin{equation}\label{claimax}
A(x)=\sum_{n_j\le x}a_j
=2x+\left(\frac12-\{(\log x)^2\}\right)\frac{x}{\log x}
+O\left(\frac{x}{(\log x)^2}\right).
\end{equation}
The Dirichlet series defining $f(s)$ converges absolutely for $\sigma>1$, the function has a continuous extension to $\sigma\ge1$, $s\ne1$, and
\[
\frac1T\int_1^T|f(1+it)|^2\,dt\sim2\log T.
\]
\end{Theorem}

The positivity of coefficients is required in Theorem \ref{thm:moment2}. To see this, take
\[
n_j=j,\qquad a_j=1+\frac{(-1)^j j}{\log^2j},\qquad j=2,3,\dots.
\]
Then
\[
A(x)=x+O\left(\frac{x}{\log^2x}\right),
\]
whereas
\[
\sum_{j=2}^{\infty}\frac{a_j^2}{j^2}
\gg\sum_{j=2}^{\infty}\frac1{\log^4j}=\infty.
\]
Moreover, techniques from \cite[Appendix B]{bh} can be adapted to show  that the mean square of the corresponding function at $\sigma=1$ tends to infinity.


The logarithmic error term is useful in some applications. Hilberdink and Lapidus \cite[Theorem 2.2]{hilberdink2006} showed that a Beurling prime number theorem with a power-saving error term yields a generalized integer-counting estimate with remainder
\begin{equation}\label{eq:introremainder}
O\left(x\exp\left(-c\sqrt{\log x\log\log x}\right)\right),\qquad c>0.
\end{equation}
This bound does not imply $O(x^\beta)$ for any fixed $\beta<1$, so the mean-square theorem of \cite{bh} does not apply directly. However, the bound \eqref{eq:introremainder} does imply \eqref{eq:Rxbound} for every fixed $\varepsilon>0$. Theorem \ref{thm:moment2} therefore shows that the second moment of the Beurling zeta function attached to such a Beurling system exists at $\sigma=1$, see Corollary \ref{thm:moment1} below. Broucke, Debruyne, and Vindas \cite{bdv20} constructed Beurling systems showing that a remainder of the size \eqref{eq:introremainder} cannot in general be replaced by a fixed power saving. For earlier investigations of the second moment of Beurling zeta functions see \cite{hilberdink2010}, \cite{dgn2019}, and \cite[Lemma 8]{drungilas2021}.

Our main application concerns the Selberg zeta function. Let $\Gamma\subset\operatorname{PSL}_2(\R)$ be any cofinite Fuchsian group. For $\sigma>1$, its Selberg zeta function is defined by the absolutely convergent Euler product
\[
Z(s)=\prod_{\{P\}}\prod_{k=0}^{\infty}(1-N(P)^{-s-k}),
\]
where $\{P\}$ runs through the primitive hyperbolic conjugacy classes of $\Gamma$ and $N(P)=\lambda^2$ if the eigenvalues of $P$ are $\lambda$ and $\lambda^{-1}$, with $|\lambda|>1$. The function $Z(s)$ has a meromorphic continuation to the whole complex plane. For more information, see Iwaniec \cite[Section 10.8]{iwaniec1995} and Jorgenson and Smajlovi\'c \cite{js2017}.

For $\sigma>1$, we can write $Z(s)$ and its multiplicative inverse as general Dirichlet series,
\begin{equation}\label{eq:selbergseries}
Z(s)=\sum_{n=1}^{\infty}\frac{b_n}{y_n^s},\qquad
\frac1{Z(s)}=\sum_{n=1}^{\infty}\frac{c_n}{y_n^s},
\end{equation}
where $1=y_1<y_2<\cdots$, $b_n,c_n\in\R$, and $c_n>0$.
Expanding the Euler products gives
\begin{equation}\label{eq:ineqZcoefs}
|b_n|\le c_n
\end{equation}
for all $n \in \N.$

\begin{Corollary}\label{thm:moments}
Let $Z(s)$ be the Selberg zeta function attached to any cofinite Fuchsian group. Then
\begin{equation}\label{eq:momentselberg}
\lim_{T\to\infty}\frac1T\int_1^T|Z(1+it)|^2\,dt
=\sum_{n=1}^{\infty}\frac{b_n^2}{y_n^2}
\end{equation}
and
\begin{equation}\label{eq:momentselberginverse}
\lim_{T\to\infty}\frac1T\int_1^T|Z(1+it)|^{-2}\,dt
=\sum_{n=1}^{\infty}\frac{c_n^2}{y_n^2}.
\end{equation}
Both series on the right-hand side converge.
\end{Corollary}

The first author, together with Drungilas and Novikas, previously investigated these moments for $\operatorname{PSL}_2(\Z)$ in \cite{drungilas2021}. Theorem 3 of that paper establishes the existence of these moments under additional assumptions, including a separation condition. Corollary \ref{thm:moments} removes these assumptions and applies to every cofinite Fuchsian group. The second moments of the logarithm and the logarithmic derivative of Selberg zeta functions were investigated in \cite{dgk2013,aoki2020,hashimoto2023,hashimoto2024}.

We establish Theorem \ref{thm:moment2} in Section \ref{proofthm:moment2}. Section \ref{sec:example} proves Theorem \ref{claim}. The applications to Beurling and Selberg zeta functions are treated in Sections \ref{BeurlingZ} and \ref{sec:momentsproof}, respectively.

We use $f(x)\ll g(x)$ or $f(x)=O(g(x))$ to mean that $|f(x)|\le C|g(x)|$ for all sufficiently large $x$, with a constant $C>0$. We write $f(x)=o(g(x))$ if $f(x)/g(x)\to0$, and $f(x)\sim g(x)$ if $f(x)/g(x)\to1$.

\section{Proof of Theorem \ref{thm:moment2}}\label{proofthm:moment2}

For the proof of Theorem \ref{thm:moment2} we adapt the ideas of Broucke and Hilberdink \cite{bh}. Let
\[
    f_N(s)=\sum_{n_j\le N}\frac{a_j}{n_j^s}.
\]
All implied constants in this section may depend on the series and on $\varepsilon$, but not on $N$ or $T$.

\begin{Lemma}\label{thm:prop1}
Assume that the conditions of Theorem \ref{thm:moment2} are satisfied. Then
\[
\frac1T\int_1^T|f_N(1+it)|^2dt
=\frac{T-1}{T}\sum_{n_j\le N}\frac{a_j^2}{n_j^2}
+O\left(\frac{\log N\log\log N}{T}\right)+o(1),
\]
for $T > 1$ and $N > e^e.$ 
Also, the $o(1)$ term is uniform in $N$.
\end{Lemma}

\begin{proof}
Writing $s(u)=\sin u/u$, we compute
\begin{equation}\label{eq:fNcomp}
\begin{split}
&\frac1T\int_1^T|f_N(1+it)|^2dt
={}\\
&\frac{T-1}{T}\sum_{n_j\le N}\frac{a_j^2}{n_j^2}
+2\sum_{n_j<n_k\le N}\frac{a_ja_k}{n_jn_k}
\left(s(T\log(n_k/n_j))
-\frac{s(\log(n_k/n_j))}{T}\right).
\end{split}
\end{equation}

Without loss of generality, we further assume that $n_j>1$. Let
$$
R_j = \frac{n_j}{(\log n_j)^{1+\varepsilon}}.
$$
Next, we consider the first double sum in the last line of \eqref{eq:fNcomp}.  
\begin{equation}\label{eq:firstdouble}
\begin{split}
&\sum_{n_j < n_k \le N} \frac{a_ja_k}{n_jn_k} s(T \log(n_k / n_j)) 
= \\
&\sum_{n_j \le N}\frac{a_j}{n_j}  
\sum_{2n_j < n_k \le N} \frac{a_k}{n_k} s(T \log(n_k / n_j)) 
+ \\
&\sum_{n_j \le N}\frac{a_j}{n_j} 
\sum_{n_j + R_j < n_k \le \min\{2n_j,N\}} \frac{a_k}{n_k} s(T \log(n_k / n_j)) + \\
&\sum_{n_j \le N}\frac{a_j}{n_j} \sum_{n_j < n_k \le\min\{n_j+R_j,N\}} \frac{a_k}{n_k} s(T \log(n_k / n_j)) =
S_1 + S_2 + S_3.
\end{split}
    \end{equation}

To estimate $S_1$ we use the bound $|s(x)| \le 1/x$ and Abel's summation formula to get
    \begin{equation*}
\begin{split}
&\sum_{2n_j \le n_k \le N} \frac{a_k}{n_k} s(T \log(n_k / n_j)) \ll \frac{1}{T} \sum_{2n_j \le n_k \le N} \frac{a_k}{n_k \log(n_k / n_j)} \ll \\
& \frac{1}{T} \left(\frac{A(N)}{N \log(N/n_j)} + \frac{A(n_j)}{n_j \log 2} + \int_{2 n_j}^N \frac{dx}{x \log(x / n_j)}  \right) \ll \frac{\log \log N}{T}.
\end{split}
    \end{equation*}
Hence,
\begin{equation}\label{eq:S1}
S_1 \ll \frac{\log \log N}{T}  \sum_{n_j \le N} \frac{a_j}{n_j} \ll \frac{\log \log N}{T}  \int_{n_1}^{N} \frac{dx}{x}  \ll \frac{\log N \log \log N}{T},
\end{equation}
where we again used Abel's summation formula.

For  $S_2$
we use  $|s(x)| \le 1/x$ and $\log(n_k/n_j) \ge (n_k - n_j)/n_k$ to show that
\begin{equation*}
\begin{split}
\sum_{n_j + R_j < n_k \le 2n_j} \frac{a_k}{n_k} s(T \log(n_k / n_j)) \ll 
\frac{1}{T}\sum_{n_j + R_j < n_k \le 2n_j} \frac{a_k}{n_k - n_j}.
\end{split}
\end{equation*}
If
$$
\frac{(\log n_j)^{1+\varepsilon}}{2} > 1,
$$
then there exists integer $K_j$ such that
$$
    \frac{n_j}{2} < 2^{K_j}R_j \le n_j.
$$
Thus
$$
K_j = \frac{1+\varepsilon}{\log 2} \log \log n_j + O(1) \ll \log \log n_j
$$
and
\begin{equation*}
\begin{split}
&\sum_{n_j + R_j < n_k \le 2n_j} \frac{a_k}{n_k - n_j} = \\
&\sum_{r=1}^{K_j} \sum_{n_j+2^{r-1}R_j < n_k\le n_j+2^{r}R_j} \frac{a_k}{n_k - n_j} + 
\sum_{n_j+2^{K_j}R_j < n_k \le 2n_j} \frac{a_k}{n_k - n_j}
\le \\
&\sum_{r=1}^{K_j} \frac{A(n_j+2^{r} R_j) - A(n_j+2^{r-1} R_j)}{2^{r-1} R_j} + 2\frac{A(2n_j) - A(n_j)}{n_j} \ll \\
&\sum_{r=1}^{K_j} \frac{2^{r} R_j}{2^{r-1} R_j} + 1 \ll K_j+1 \ll \log\log n_j.
\end{split}
\end{equation*}
This leads to
\begin{equation}\label{eq:S2}
S_2 \ll \frac{\log \log N}{T}  \sum_{n_j \le N} \frac{a_j}{n_j} \ll \frac{\log \log N}{T}  \int_{n_1}^{N} \frac{dx}{x}  \ll \frac{\log N \log \log N}{T}.
\end{equation}

We will show that $S_3 = o(1)$ uniformly in $N>e^e$.
By $|s(x)| \le 1$ we get
$$
\sum_{n_j < n_k \le n_j + R_j} \frac{a_k}{n_k} \le \frac{A(n_j+ R_j)-A(n_j)}{n_j} \ll \frac{R_j}{n_j} = \frac{1}{(\log n_j)^{1+\varepsilon}}.
$$
Thus, independent of $N$ and $T$, we have
$$
|S_3| \ll \sum_{n_j \le \infty} \frac{a_j}{n_j(\log n_j)^{1+\varepsilon}} \ll \int_{n_1}^\infty \frac{dx}{x (\log x)^{1+\varepsilon}}  \ll 1.
$$
Hence, for any $\eta > 0$, there is a positive $n_0=n_0(\eta)$  such that
$$
\left| \sum_{n_0 \le n_j \le \infty} \frac{a_j}{n_j}\sum_{n_j < n_k \le n_j + R_j} \frac{a_k}{n_k} s(T \log(n_k / n_j)) \right| \le \frac{\eta}{2}
$$
for all $T>0$. In view of 
$$
    \left| \sum_{n_j \le n_0} \sum_{n_j < n_k \le n_j + R_j} \frac{a_ja_k}{n_jn_k} s(T \log(n_k / n_j)) \right| \ll_{n_0} \frac{1}{T},
$$
there is a positive $T_0=T_0(n_0(\eta))$ such that, for all $T > T_0$,
$$
    \left| \sum_{n_j \le n_0} \sum_{n_j < n_k \le n_j + R_j} \frac{a_ja_k}{n_jn_k} s(T \log(n_k / n_j)) \right| \le \frac{\eta}{2}.
$$
Hence, if $\eta > 0$, then there is $T_0(\eta)$ such that, for $T > T_0(\eta)$ , the bound
$
|S_3| \le \eta
$ is valid,
 i.e. 
\begin{equation}\label{eq:S3}
 S_3 = o(1)
 \end{equation}
 uniformly in $N$.

Finally, split the sum containing $s(\log(n_k/n_j))/T$ in \eqref{eq:fNcomp} into the same three ranges. These ranges can be estimated in the same way as $S_1,$ $S_2$ and $S_3$ and give the same error bounds.
\end{proof}

\begin{Lemma}\label{thm:prop2}
Assume that the conditions of Theorem \ref{thm:moment2} are satisfied. Then
\[
\frac1T\int_1^T|f(1+it)-f_N(1+it)|^2dt
\ll\frac1T+\frac1{(\log N)^{2+2\varepsilon}}
+\frac{T}{(\log N)^{1+2\varepsilon}},
\]
when $N>e^e$ and $T>1$.
\end{Lemma}

\begin{proof}
Since $A(x)\ll x$, the Dirichlet series for $f(s)$ converges absolutely for $\sigma>1$. Abel's summation formula gives
\begin{equation}\label{eq:diffexp}
f(s)-f_N(s)
=\frac{\varrho N^{1-s}}{s-1}-\frac{R(N)}{N^s}
+s\int_N^\infty\frac{R(x)}{x^{s+1}}dx.
\end{equation}
The integral converges absolutely for $\sigma\ge1$, because
\[
\int_N^\infty\frac{R(x)}{x^{s+1}}dx \ll \int_N^\infty\frac{dx}{x(\log x)^{1+\varepsilon}} \ll 1.
\]
Thus \eqref{eq:diffexp} gives the continuous extension of $f(s)$ to $\sigma\ge1$, $s\ne1$. For $1\le t\le T$, it follows that
\begin{equation}\label{fn1}
|f(1+it)-f_N(1+it)|
\ll\frac1t+\frac1{(\log N)^{1+\varepsilon}}
+T\left|\int_N^\infty\frac{R(x)}{x^{2+it}}dx\right|.
\end{equation}
Squaring and integrating, we obtain
\begin{equation}\label{boundwithdoubleint}
\begin{split}
&\frac1T\int_1^T|f(1+it)-f_N(1+it)|^2dt
\ll\\
&\frac1T+\frac1{(\log N)^{2+2\varepsilon}}
+T\int_1^T\left|\int_N^\infty\frac{R(x)}{x^{2+it}}dx\right|^2dt.
\end{split}
\end{equation}
Set
\[
h_N(u)=\begin{cases}
    e^{-u}R(e^u), \text{ if } u \ge \log N, \\
    0, \text{ otherwise}.
\end{cases}
\]
Making the change of variable $x=e^u$ we get
\[
\int_N^\infty\frac{R(x)}{x^{2+it}}dx
=\int_{- \infty}^\infty h_N(u)e^{-itu}du.
\]
That is, as a function of $t,$ the integral
$$
\int_N^\infty\frac{R(x)}{x^{2+it}}dx
$$
is the Fourier transform of $h_N(u).$
The function $h_N(u)$ belongs to $L^1(\R)\cap L^2(\R)$, so we can apply Plancherel's theorem to get
\begin{equation}\label{eq:plancherel}
\begin{split}
&\int_1^T\left|\int_N^\infty\frac{R(x)}{x^{2+it}}dx\right|^2dt
\ll \int_{\log N}^\infty|e^{-u}R(e^u)|^2du\\
&\ll\int_{\log N}^\infty\frac{du}{u^{2+2\varepsilon}}
\ll\frac1{(\log N)^{1+2\varepsilon}}.
\end{split}
\end{equation}
The lemma then follows from \eqref{boundwithdoubleint}.
\end{proof}

\begin{proof}[Proof of Theorem \ref{thm:moment2}]
Absolute convergence for $\sigma>1$ and the continuous extension to $\sigma\ge1$, $s\ne1$, were shown in the proof of Lemma \ref{thm:prop2}.

The series on the right-hand side of \eqref{mainmoment} converges because
    $$
\sum_{j=1}^{\infty}\frac{a^2_j}{n_j^2} \ll\sum_{j=1}^{\infty}\frac{a_j(A(n_j)-A(n_{j}-1))}{n_j^2} \ll
\sum_{j=1}^{\infty}\frac{a_j}{n_j(\log n_j)^{1+\varepsilon}}.
    $$ 

Next, we rewrite the LHS of \eqref{mainmoment} as
\begin{equation}\label{eq:polarization}
\begin{split}
&\frac{1}{T} \int_1^T\left|f(1+it)\right|^2dt = \frac{1}{T} \int_1^T|f_N(1+it)|^2dt +
\\&  \frac{1}{T}\int_1^T\left|f(1+it) - f_N(1+it)\right|^2dt + \\
&2\Re \left(\frac{1}{T}\int_1^T f_N(1+it)\left(\overline{f(1+it)-f_N(1+it)}\right)dt\right),
\end{split}
\end{equation}
where, by the Cauchy-Schwarz inequality, the last term is bounded by
    $$
        2\left(\frac{1}{T^2}\int_1^T|f_N(1+it)|^2dt\int_1^T|f(1+it) - f_N(1+it)|^2dt\right)^{1/2}.
    $$
Choose $N=N(T)$ such that 
$
 T \sim (\log N)^{1+\varepsilon}.
$
Then using Lemmas \ref{thm:prop1} and \ref{thm:prop2} we obtain the moment \eqref{mainmoment}.
This concludes the proof.
\end{proof}

The following corollary allows us to treat general Dirichlet series with real coefficients which are not necessarily positive.

\begin{Corollary}\label{cor:signed}
Suppose that $f(s)=\sum_j a_jn_j^{-s}$ satisfies the assumptions of Theorem \ref{thm:moment2}. Let $a'_j\in\R$, $|a'_j|\le a_j$, and assume that, for some $\delta>0$,
\[
B(x)=\sum_{n_j\le x}a'_j
=O\left(\frac{x}{(\log x)^{1+\delta}}\right).
\]
Then the series
\[
g(s)=\sum_{j=1}^\infty\frac{a'_j}{n_j^s}
\]
converges absolutely for $\sigma>1$ and has a continuous extension to $\sigma\ge1$, $s\ne1$. Moreover,
\[
\lim_{T\to\infty}\frac1T\int_1^T|g(1+it)|^2dt
=\sum_{j=1}^\infty\frac{(a'_j)^2}{n_j^2},
\]
and the series on the right-hand side converges.
\end{Corollary}

\begin{proof}
Absolute convergence and convergence of 
$$
\sum_{j=1}^\infty\frac{(a'_j)^2}{n_j^2}
$$
follow by comparison with $f$. Abel's summation formula gives
\[
g(s)-\sum_{n_j\le N}\frac{a'_j}{n_j^s}
=-\frac{B(N)}{N^s}+s\int_N^\infty\frac{B(x)}{x^{s+1}}dx,
\]
whose integral converges absolutely for $\sigma\ge1$. This proves the asserted continuous extension.

The coefficients of $f+g$ and $f-g$ are nonnegative, and their summatory functions satisfy
\[
\sum_{n_j\le x}(a_j\pm a'_j)
=\varrho x+O\left(\frac{x}{(\log x)^{1+\min\{\varepsilon,\delta\}}}\right).
\]
After omitting zero coefficients, Theorem \ref{thm:moment2} applies to both of these series.
The identity
\[
|f+g|^2+|f-g|^2=2|f|^2+2|g|^2
\]
and the three mean-square formulas for $f+g$, $f-g$, and $f$ now give
\begin{align*}
&\lim_{T\to\infty}\frac1T\int_1^T|g(1+it)|^2dt = \\
&\frac12\sum_j\frac{(a_j+a'_j)^2+(a_j-a'_j)^2}{n_j^2}
-\sum_j\frac{a_j^2}{n_j^2}
=\sum_j\frac{(a'_j)^2}{n_j^2}.
\end{align*}
\end{proof}

\section{Proof of Theorem \ref{claim}}\label{sec:example}

In this section, we use the notation
\[
n_j=e^{\sqrt j},\qquad a_j=\frac{e^{\sqrt j}}{\sqrt j},\qquad
A(x)=\sum_{n_j\le x}a_j,
\]
and
\[
f(s)=\sum_{j=1}^{\infty}\frac{e^{-(s-1)\sqrt j}}{\sqrt j}
\quad(\sigma>1),\qquad
f_N(s)=\sum_{n_j\le N}\frac{a_j}{n_j^s}.
\]
We first establish the continuous extension and the estimate for $A(x)$. We then compare $f$ with $f_N$ in mean square.

\Needspace{6\baselineskip}
\begin{Lemma}\label{claimextension}
The function $f(s)$ has a continuous extension to $\sigma\ge1$, $s\ne1$.
\end{Lemma}

\begin{proof}
Change variable $w=s-1$ and set
\[
g(w)=\sum_{j=1}^{\infty}\frac{e^{-w\sqrt j}}{\sqrt j},\qquad \Re w>0.
\]
Applying the Euler--Maclaurin formula to the partial sums and then letting the upper endpoint tend to infinity gives
\begin{equation}\label{eq:examplesumRHS}
g(w)=\frac{2e^{-w}}{w}+\frac{e^{-w}}2
+\frac{e^{-w}(w+1)}{24}
-\frac12\int_1^{\infty}B(x)h_w(x)\,dx,
\end{equation}
where
\[
B(x)=\{x\}^2-\{x\}+\frac16,\qquad
h_w(x)=\frac{e^{-w\sqrt x}(w^2x+3w\sqrt x+3)}{4x^{5/2}}.
\]
Here $h_w$ is the second derivative of $x^{-1/2}e^{-w\sqrt x}$, and
$\int_1^{\infty}x^{-1/2}e^{-w\sqrt x}\,dx=2e^{-w}/w$ for $\Re w>0$.
For any compact subset $K$ of $\Re w\ge0$,
\[
B(x)h_w(x)\ll_K x^{-3/2},
\]
for $w \in K$.
Thus, the dominated convergence theorem shows that the integral in \eqref{eq:examplesumRHS} is continuous up to $\Re w=0$. The right-hand side gives the required extension when $w\ne0$.
\end{proof}

\begin{Lemma}\label{thm:Axestimate}
We have
\[
A(x)=2x+\left(\frac12-\{(\log x)^2\}\right)\frac{x}{\log x}
+O\left(\frac{x}{(\log x)^2}\right).
\]
\end{Lemma}

\begin{proof}
Let $\floor{u} = u - \{ u \}.$ Abel's summation formula gives
\begin{align*}
&A(x)=\sum_{j\le (\log x)^2}\frac{e^{\sqrt j}}{\sqrt j}\\
&=\floor{(\log x)^2}\frac{x}{\log x}
-\int_1^{(\log x)^2}\frac{\floor{u}e^{\sqrt u}}{2u}\,du
+\int_1^{(\log x)^2}\frac{\floor{u}e^{\sqrt u}}{2u^{3/2}}\,du\\
&=2x-\{(\log x)^2\}\frac{x}{\log x}
+\int_1^{(\log x)^2}\frac{\{u\}e^{\sqrt u}}{2u}\,du
-\int_1^{(\log x)^2}\frac{\{u\}e^{\sqrt u}}{2u^{3/2}}\,du+O(1).
\end{align*}
In the last step we used
\[
\int_1^{(\log x)^2}\frac{e^{\sqrt u}}2\,du=x\log x-x
\text{ and }
\int_1^{(\log x)^2}\frac{e^{\sqrt u}}{2\sqrt u}\,du=x+O(1).
\]
Write $\{u\}=1/2+(\{u\}-1/2)$ and integrate by parts to get
\begin{align*}
&A(x)=2x+\left(\frac12-\{(\log x)^2\}\right)\frac{x}{\log x}\\
&+\int_1^{(\log x)^2}\frac{(\{u\}-1/2)e^{\sqrt u}}{2u}\,du
-\int_1^{(\log x)^2}\frac{(\{u\}-1/2)e^{\sqrt u}}{2u^{3/2}}\,du+O(1).
\end{align*}
The function
\[
\varphi(v)=\int_0^v(\{u\}-1/2)\,du
\]
is bounded. Another integration by parts therefore gives, for $q=1$ and $q=3/2$,
\begin{align*}
&\int_1^{(\log x)^2}(\{u\}-1/2)\frac{e^{\sqrt u}}{2u^q}\,du=\\
&\varphi((\log x)^2)\frac{x}{2(\log x)^{2q}} - \varphi(1)\frac{e}{2}
-\int_1^{(\log x)^2}\varphi(u)\,d\left(\frac{e^{\sqrt u}}{2u^q}\right)
\ll\frac{x}{(\log x)^2}.
\end{align*}
\end{proof}

\begin{Lemma}\label{thm:exampletruncationestimate}
Uniformly for $N>e^e$ and $T>1$,
\[
\frac1T\int_1^T|f_N(1+it)|^2\,dt
=2\frac{T-1}{T}\log\log N+O(1)
+O\left(\frac{\log N\log\log N}{T}\right).
\]
\end{Lemma}

\begin{proof}
As in \eqref{eq:fNcomp}, write
\begin{equation*}
\begin{split}
&\frac1T\int_1^T|f_N(1+it)|^2dt
={}\\
&\frac{T-1}{T}\sum_{n_j\le N}\frac{a_j^2}{n_j^2}
+2\sum_{n_j<n_k\le N}\frac{a_ja_k}{n_jn_k}
\left(s(T\log(n_k/n_j))
-\frac{s(\log(n_k/n_j))}{T}\right),
\end{split}
\end{equation*}
where $s(x)=\sin(x)/x.$
We have
$$
\sum_{n_j \le N} \frac{a_j^2}{n_j^2} = \sum_{j \le (\log N)^2} \frac{1}{j} = 2 \log \log N + O(1).
$$
Using $|s(x)| \le 1/x$, it remains to prove that
$$
    \sum_{n_j < n_k \le N} \frac{a_j a_k}{n_j n_k \log(n_k / n_j)} \ll \log N \log \log N.
$$
This sum is equal to 
$$
    \sum_{j < k \le ( \log N)^2} \frac{1}{\sqrt{j} \sqrt{k} (\sqrt{k} - \sqrt{j})}.
$$
We see that
$$
\sum_{k=j+1}^{(\log N)^2} \frac{1}{\sqrt{k} (\sqrt{k} - \sqrt{j})} = \sum_{k=j+1}^{(\log N)^2} \frac{\sqrt{k}+\sqrt{j}}{\sqrt{k}(k-j)} \ll \log \log N.
$$
Hence
$$
\sum_{j < k \le ( \log N)^2} \frac{1}{\sqrt{j} \sqrt{k} (\sqrt{k} - \sqrt{j})} \ll \log \log N \sum_{j \le ( \log N)^2}\frac{1}{\sqrt{j}} \ll \log N \log \log N.
$$
\end{proof}

\begin{Lemma}\label{thm:examplediffbound}
Uniformly for $N>e^e$ and $T>1$,
\[
\frac1T\int_1^T|f(1+it)-f_N(1+it)|^2\,dt
\ll\frac1T+\frac1{(\log N)^2}+\frac{T}{\log N}.
\]
\end{Lemma}

\begin{proof}
By Lemma \ref{thm:Axestimate}, write $A(x)=2x+R(x)$, where $R(x)\ll x/\log x$. For $1<\sigma\le2$, Abel's summation formula gives
\[
f(s)-f_N(s)=\frac{2N^{1-s}}{s-1}-\frac{R(N)}{N^s}
+s\int_N^{\infty}\frac{R(x)}{x^{s+1}}\,dx.
\]
As in \eqref{boundwithdoubleint}, for $T>1$ we obtain
\begin{align*}
&\frac1T\int_1^T|f(\sigma+it)-f_N(\sigma+it)|^2\,dt
\ll\\
&\frac1T+\frac1{(\log N)^2}
+T\int_1^T\left|\int_N^{\infty}\frac{R(x)}{x^{\sigma+1+it}}\,dx\right|^2dt.
\end{align*}
The change of variable $x=e^u$ and Plancherel's theorem give (see \eqref{eq:plancherel})
\begin{align*}
&\int_1^T\left|\int_N^{\infty}\frac{R(x)}{x^{\sigma+1+it}}\,dx\right|^2dt
\ll\pi\int_{\log N}^{\infty}|e^{-\sigma u}R(e^u)|^2\,du\\
&\ll\int_{\log N}^{\infty}\frac{e^{-2(\sigma-1)u}}{u^2}\,du
\ll\frac{1}{\log N}.
\end{align*}
All constants are independent of $1<\sigma\le2$. Lemma \ref{claimextension} gives the continuous extension of $f$ to $\sigma=1$. Letting $\sigma\rightarrow1$ from the right and applying Fatou's lemma proves the result.
\end{proof}

\begin{proof}[Proof of Theorem \ref{claim}]
The continuous extension is obtained in Lemma \ref{claimextension}, and \eqref{claimax} follows from Lemma \ref{thm:Axestimate}.

Let $N = N(T)$ be a sequence such that
$
    T \sim \log N (\log \log N)^{1/2}.
$
Then, similarly as in \eqref{eq:polarization},
\begin{equation*}
    \begin{split}
&\frac{1}{T} \int_1^T\left|f(1+it)\right|^2dt= \\
&\frac{1}{T} \int_1^T|f_N(1+it)|^2dt +  \frac{1}{T}\int_1^T\left|f(1+it) - f_N(1+it)\right|^2dt
\\&+  
O\left(\left(\frac{1}{T^2}\int_1^T|f_N(1+it)|^2dt\int_1^T|f(1+it) - f_N(1+it)|^2dt\right)^{1/2}\right).
\end{split}
\end{equation*}
By Lemma \ref{thm:exampletruncationestimate} the first integral on the right hand side is
$$
 2\frac{T-1}{T} \log \log N + O\left((\log \log N)^{1/2}\right),
$$
by Lemma \ref{thm:examplediffbound} the second integral is
$
O\left((\log \log N)^{1/2}\right),
$
and the error term is
$
O\left((\log \log N)^{3/4}\right).
$
Hence,
$$
\frac{1}{T} \int_1^T\left|f(1+it)\right|^2dt= 2\frac{T-1}{T}\log \log N + O\left((\log \log N)^{3/4}\right).
$$
By our choice of $N$ we get
$$
\log T \sim \log \log N + \frac{1}{2} \log \log \log N.
$$
This proves the claim.
\end{proof}

\section{Beurling zeta functions}\label{BeurlingZ}

In this section we consider a generalized prime number system $\mathcal P$, often called $g$-primes,
\[
1<p_1\le p_2\le\dots,
\]
with $p_j\in\mathbb R$ tending to infinity, and the associated Beurling zeta function
\begin{equation}\label{eulerbeurling}
\zeta_{\mathcal P}(s)=\prod_{j=1}^{\infty}\left(1-\frac1{p_j^s}\right)^{-1}.
\end{equation}
A formal multiplication of the product gives the Dirichlet series
\begin{equation}\label{dirichletbeurling}
\zeta_{\mathcal P}(s)=\sum_{n=1}^\infty\frac{d_n}{\nu_n^s},
\end{equation}
where $1=\nu_1<\nu_2<\dots$ is the increasing sequence of power-products of $g$-primes, with corresponding multiplicities $1=d_1,d_2,\dots$.
We assume that, for some fixed $1/2<\alpha<1$,
\begin{equation}\label{pxli}
\pi_{\mathcal P}(x)=\sum_{p_j\le x}1=\operatorname{Li}(x)+O(x^\alpha).
\end{equation}
In particular, both \eqref{eulerbeurling} and \eqref{dirichletbeurling} converge absolutely for $\sigma>1$.

To apply the results of Hilberdink and Lapidus, we first pass to the weighted prime-power counting function
\[
\psi_{\mathcal P}(x)=\sum_{p_j^m\le x,\ m\ge1}\log p_j.
\]
Abel's summation formula gives
\[
\sum_{p_j\le x}\log p_j=x+O(x^\alpha\log x).
\]
The terms with $m\ge2$ contribute $O(x^{1/2}\log x)$. Hence, for every fixed $\theta\in(\alpha,1)$,
\begin{equation}\label{eq:beurlingpsi}
\psi_{\mathcal P}(x)=x+O(x^\alpha\log x)=x+O(x^\theta).
\end{equation}
The first estimate also gives $\psi_{\mathcal P}(x)=x+O(x^{\alpha+\delta})$ for every $\delta>0$. Thus \cite[Theorem~2.1]{hilberdink2006} shows that $\zeta_{\mathcal P}$ has a meromorphic continuation to $\sigma>\alpha$, with a simple pole at $s=1$ and no zeros in this half-plane. Applying \cite[Theorem~2.2]{hilberdink2006} with the fixed exponent $\theta$ in \eqref{eq:beurlingpsi}, we obtain
\begin{equation}\label{betax}
\sum_{\nu_n\le x}d_n=\eta x+D(x),
\end{equation}
where $\eta>0$ and, for some fixed $c>0$,
\begin{equation}\label{Dx}
D(x)\ll x\exp\left(-c\sqrt{\log x\log\log x}\right).
\end{equation}
This remainder satisfies \eqref{eq:Rxbound} for every fixed $\varepsilon>0$. Consequently, Theorem~\ref{thm:moment2} gives the following result.

\begin{Corollary}\label{thm:moment1}
Under assumption \eqref{pxli}, we have
\begin{equation}\label{eq:momentzeta1}
\lim_{T\to\infty}\frac1T\int_1^T|\zeta_{\mathcal P}(1+it)|^2\,dt
=\sum_{n=1}^{\infty}\frac{d_n^2}{\nu_n^2}.
\end{equation}
The series on the right-hand side converges.
\end{Corollary}

Note that Broucke, Debruyne, and Vindas \cite{bdv20} constructed a Beurling generalized number system which satisfies \eqref{pxli} with $\alpha=1/2$ and \eqref{betax} with
\[
D(x)=\Omega\left(x\exp\left(-C\sqrt{\log x\log\log x}\right)\right)
\]
for some $C>0$. Thus a prime number theorem of the form \eqref{pxli} does not, in general, give a power-saving remainder for the generalized integer counting function. Corollary~\ref{thm:moment1} removes the separation hypothesis in \cite[Lemma 8]{drungilas2021}, where the second moment of $\zeta_{\mathcal P}(1+it)$ was treated. 

\section{Selberg zeta functions}\label{sec:momentsproof}

Let $Z(s)$ be the Selberg zeta function for a cofinite Fuchsian group $\Gamma$, and retain the coefficients $b_n,c_n$ and the sequence $(y_n)$ from \eqref{eq:selbergseries}. In particular, $c_n>0$ and $|b_n|\le c_n$. To prove Corollary~\ref{thm:moments}, we establish the following coefficient estimates.

\begin{Lemma}\label{thm:ZetaCoefEstimateLemma}
For some fixed $c>0$, we have
\begin{equation}\label{eq:ZetaCoefEstimate}
\sum_{y_n\le x}b_n
=O\left(x\exp\left(-c\sqrt{\log x\log\log x}\right)\right)
\end{equation}
and
\begin{equation}\label{eq:Zeta-1CoefEstimate}
\sum_{y_n\le x}c_n
=\varrho x+O\left(x\exp\left(-c\sqrt{\log x\log\log x}\right)\right),
\end{equation}
where $\varrho>0$ is the residue of $1/Z(s)$ at $s=1$.
\end{Lemma}

\begin{proof}
For $\sigma > 1,$ we can write
$$
    Z(s) = Z(s+1)\prod_{P}(1-N(P)^{-s}).
$$
Corollary~\ref{thm:moments} will follow by considering the Ruelle zeta function, which is defined as
\begin{align}\label{ruelle}
Z_1(s)=Z(s)/Z(s+1)=\prod_{P}(1-N(P)^{-s}).
\end{align}
This function has a meromorphic continuation to $\C$.

Let $Z_2(s) = Z(s+1)$.
Since $Z(s)$ is defined for $\sigma > 1$ by an absolutely convergent Euler product, there exists a constant $C > 0$ such that
\begin{align}\label{Z2bound}
        |Z_2(s)|^{\pm1} \le C,
\end{align}
for all $\sigma > 1/2.$

Note that for $N(P)$, the prime geodesic theorem is valid, i.e. there exists $0 < \alpha < 1$ such that 
\begin{equation}\label{eq:primegeodesictheorem}
    \sum_{N(P) \le x} 1 = \operatorname{Li}(x)+ O(x^\alpha).
\end{equation}
A suitable form of the prime geodesic theorem can be found in Iwaniec \cite[Theorem 10.5]{iwaniec1995}.

The function $1/Z_1(s)$ is a Beurling zeta-function (see \cite[Ending notes]{dgn2019} and \cite[Introduction]{drungilas2021}).
In view of \eqref{eq:primegeodesictheorem}, we can adapt the proof of \cite[Theorem 2.2]{hilberdink2006}.

In \cite[Proof of Theorem 2.2, between formulas (2.5) and (2.6)]{hilberdink2006} we find that, for $\sigma\in(\alpha, 1)$,
$$
    \log 1/Z_1(\sigma+it) = -\int_{\sigma}^2 \frac{(1/Z_1)'(y+it)}{1/Z_1(y+it)}dy + O(1)
$$
and, for $|t|\ge 3$ and $\varepsilon(t) = (1-\alpha)\frac{\log \log|t|}{\log |t|}$, 
$$
\int_{1-\varepsilon(t)}^2 \left|\frac{(1/Z_1)'(y+it)}{1/Z_1(y+it)} \right|dy \le (1 + o(1)) \frac{\log|t|}{\log\log|t|}. 
$$
Consequently,
\begin{align}\label{pmBound}
    |Z_1\left(1-\varepsilon(t)+it\right)|^{\pm1} \le \exp\left(\frac{2\log|t|}{\log \log |t|}\right)
\end{align}
for all $|t|$ sufficiently large.

Denote
\begin{align*}
S(x) =\sum_{y_n \le x}c_n.
\end{align*}
Integrating Perron's formula gives
$$
    S_1(x)=\int_0^x S(y)dy = \frac{1}{2\pi i} \int_{c-i\infty}^{c+i\infty} \frac{1/Z(s)}{s(s+1)}x^{s+1}ds \qquad (c>1).
$$
Shifting the contour past the simple pole at $s=1$ of $1/Z(s)$ to the left to the contour $\gamma$, defined by $s = 1 - \varepsilon(t) + it$, when $|t| \ge 3$, and $s = 1-\varepsilon(3) + it,$ when $|t| < 3$, we get
$$
    S_1(x) = \frac{\varrho}{2}x^2 + \frac{1}{2\pi i} \int_{\gamma} \frac{1/Z(s)}{s(s+1)}x^{s+1}ds.
$$
Then the bounds \eqref{Z2bound} and \eqref{pmBound} lead to
\begin{align*}
&\left|S_1(x) - \frac{\varrho}{2}x^2 \right| \ll \\
&\int_3^\infty \frac{|1/Z_1\left(1-\varepsilon(t)+it\right)||1/Z_2\left(1-\varepsilon(t)+it\right)|}{t^2}x^{2-\varepsilon(t)}dt+ x^{2-\varepsilon(3)} \ll \\
    &x^2\int_3^\infty\exp \left(\frac{2\log t}{\log \log t} - 2\log t - (1-\alpha)\frac{\log\log t \log x}{\log t}\right) dt+ x^{2-\varepsilon(3)}. 
\end{align*}
The last integral is estimated in the proof of \cite[Theorem 2.2, see (2.7) and below]{hilberdink2006} to be
$$
    O\left(\exp\left(-d\sqrt{\log x \log \log x}\right)\right),
$$
for some constant $d > 0$.
Thus, we arrive at
$$
    S_1(x) = \frac{\varrho}{2}x^2 + O\left(x^2\exp\left(-d\sqrt{\log x \log \log x}\right)\right).
$$
Let
 $$y = x\exp\left(-\frac{d}2\sqrt{\log x \log \log x}\right).$$
Since $S(x)$ is non-decreasing, we have
\begin{align}\label{nondecrsingS}
\frac{S_1(x)-S_1(x-y)}{y}\le S(x)\le  \frac{S_1(x+y)-S_1(x)}{y},  
\end{align}
which yields the bound \eqref{eq:Zeta-1CoefEstimate}.

It is not clear how to evaluate \eqref{eq:ZetaCoefEstimate} directly, because coefficients $b_n$ may be negative. Instead, we consider 
\begin{align*}
\sum_{y_n \le x}(b_n+c_n),
\end{align*}
where $b_n+c_n\ge0$, see \eqref{eq:ineqZcoefs}. This non-negativity allows us to apply the argument used in \eqref{nondecrsingS}. In view of \eqref{eq:Zeta-1CoefEstimate}, to show \eqref{eq:ZetaCoefEstimate} it is enough to prove that
\begin{align*}
 \sum_{y_n \le x}(b_n+c_n)=\varrho x +   O\left(x\exp\left(-c\sqrt{\log x \log\log x}\right)\right). 
\end{align*}

We can do this using the same method we applied previously.
We get
$$
    \int_{0}^x \sum_{y_n \le y}(b_n+c_n) dy  = \frac{\varrho}{2}x^2 + \frac{1}{2\pi i} \int_{\gamma} \frac{Z(s) + 1/Z(s)}{s(s+1)}x^{s+1}ds.
$$
Estimating the right-hand side applying the  bounds \eqref{Z2bound} and \eqref{pmBound} to both $Z$ and $1/Z$, we get the desired asymptotics.
\end{proof}

\begin{proof}[Proof of Corollary~\ref{thm:moments}]
The estimate \eqref{eq:Zeta-1CoefEstimate} satisfies the hypothesis of Theorem~\ref{thm:moment2} for every fixed $\varepsilon>0$. Applying that theorem to $1/Z(s)$ proves \eqref{eq:momentselberginverse} and the convergence of its coefficient-square series. The estimate \eqref{eq:ZetaCoefEstimate}, together with $|b_n|\le c_n$, allows us to apply Corollary~\ref{cor:signed} to $Z(s)$. The corollary follows.
\end{proof}

Acknowledgment. This work is funded by the Research Council of Lithuania (LMTLT), agreement No. S-MIP-22-81.

\bibliographystyle{plain}
\bibliography{Zetamoment}

\end{document}